\documentclass[12pt,a4paper]{article}
\usepackage{amsmath,latexsym,amsfonts,amssymb}

\def\bql{\begin{equation}\label}
\def\eql{\end{equation}\noindent}

\def\bew{\par\noindent{\bf Proof }}
\def\qed{\quad\hfill\mbox{\P}}
\def\qend{\hfill\mbox{$\lozenge$}}

\def\p{\pi}
\def\q{\theta}
\def\s{\sigma}

\def\FC{{\cal F}}
\def\WG{\Omega}

\def\pbb{{\mathbb P}}

\def\nf{\infty}

\def\sm{\setminus}
\def\ss{\subset}
\def\top{\buildrel\pbb\over\rightarrow}

\def\noi{\noindent}

\newtheorem{theorem}{Theorem}[section]
\newtheorem{lemma}[theorem]{Lemma}
\newtheorem{proposition}[theorem]{Proposition}
\newtheorem{remark}{Remark}

\begin{document}
\date{\today}
\title{The Borel-Cantelli Lemmas for contaminated events with an application to small maxima}
\author{Guus Balkema\\
(University of Amsterdam)}
\maketitle
\setcounter{page}{1}
\begin{abstract}
For a sequence of independent events $A_n$ the sum of the associated zero-one random variables $1_{A_n}$ is almost surely finite or almost surely infinite according as the sum of the probabilities converges or diverges. In this paper the events are contaminated. What can one say about $\sum1_{B_n}$ when $B_n=A_n\sm E_n$ for a sequence of events $E_n$ with vanishing probability? It will be shown that $\sum1_{B_n}$ is infinite almost surely if $\sum PA_n=\nf$ and $E_n$ is independent of $A_n$. 
\end{abstract}

\noi{\bf Keywords:} contamination, partial maxima

\noi{\bf 2010 MSC:} 60F15 60F20 60E15

\setcounter{section}{0}
\setcounter{equation}{0}

\section{Borel-Cantelli lemmas for contaminated events}

Assume the events $A_n$ are independent and $E_n$ is independent of $A_n$. If $\pbb E_n\to0$ the probability of the events $A_n$ and $B_n=A_n\sm E_n$ are asymptotically equal. The sequence $\sum\pbb B_n$ diverges if and only if $\sum\pbb A_n$ diverges. Divergence of $\sum\pbb A_n$ implies almost sure divergence of $\sum1_{A_n}$ by the second Lemma of Borel-Cantelli. It will be shown that this also implies almost sure divergence of $\sum1_{B_n}$.

\begin{theorem}\label{tt}
Let $A_n$ be independent events and let the event $E_n$ be independent of $A_n$ for each index $n$. Assume $\pbb E_n\to0$. Set $B_n=A_n\sm E_n$.

\noi1) If $\sum\pbb B_n<\nf$ then $\sum1_{A_n}<\nf$ almost surely.

\noi2) If $\sum\pbb A_n=\nf$ then $\sum1_{B_n}=\nf$ almost\ surely.
\end{theorem}

\bew The first statement is obvious by asymptotic equality. It is included for the sake of symmetry. The second statement follows from the proposition below. \qed

\begin{remark}
The second statement is sharp: Suppose $E$ is independent of the sequence $(A_n)$ and  has positive probability. Take $E_n=E$ for all $n$. Then $\sum1_{B_n}$ vanishes on $E$. Hence the condition $\pbb E_n\to0$ is necessary. 
Independence of the events $E_n$ and $A_n$ cannot be dropped either. If $E$ has positive probability and we assume that $\pbb A_n\to0$ the events $E_n=E\cap A_n$ have vanishing probability  but the conclusion in 2) is invalid since $\sum1_{B_n}$ vanishes on $E$. \qend
\end{remark}

\begin{remark}
The proposition below allows us to relax the condition
\bql{q1}
\pbb(E_n\cap A_n)=\pbb A_n\pbb E_n\qquad n\ge1
\eql
in Theorem~\ref{tt} to
\bql{q2}
\pbb(E_n\cap A_n)=O(\pbb A_n\pbb E_n)\qquad n\to\nf.
\eql
\end{remark}

\begin{remark}
In the application below there is a filtration $\FC_0\ss\FC_1\ss\ldots$ such that $A_n\in\FC_n$ is independent of $\FC_{n-1}$ (the past) for $n\ge1$, and $E_n\in\FC_{n-1}$. This suggests a proof based on Paul L\'evy's powerful extension of the second Borel-Cantelli Lemma.  A proof which uses Serfling's Theorem, see~\cite{C12}, along these lines is possible, but does not do justice to the triviality of the result. \qend
\end{remark}

\begin{proposition}
Let $A_n$ be independent events. Let $E_n$ be events and write
\bql{q3}
\pbb(E_n\cap A_n)=e_n\pbb A_n.
\eql
If $\sum\pbb A_n=\nf$ and $e_n\to0$ the events $B_n=A_n\sm E_n$ occur infinitely often almost surely.
\end{proposition}
\bew Set $p_n=\pbb A_n$ and $D_n=E_n\cap A_n$. If $\sum\pbb D_n<\nf$ then $\sum1_{A_n}=\nf$ almost surely and $\sum1_{D_n}<\nf$ almost surely implies $\sum1_{B_n}=\nf$ almost surely, as desired. So introduce a sequence of independent  variables $U_n$ uniformly distributed on $(0,1)$, and independent of the $\s$-algebra generated by the events $A_n, E_n$, $n\ge1$. (Replace the probability space $\WG$ by $\WG\times(0,1)^\nf$ if need be.) By the lemma below one may choose $q_n\in[0,1]$ such that $\sum q_np_n=\nf$ and $\sum q_np_ne_n<\nf$. The events $A'_n=A_n\cap\{U_n\le q_n\}$ are independent with probability $p'_n=q_np_n$. Hence $\sum1_{A'_n}=\nf$ almost surely. The events $D_n\cap\{U_n\le q_n\}$ have probability $q_np_ne_n$ with finite sum. It follows that the events $B_n\cap\{U_n\le q_n\}$ almost surely occur infinitely often, and hence so do the events $B_n$. \qed

\begin{lemma}
Let $p_n$ and $a_n$ be non-negative, $\sum p_n=\nf$ and $a_n\to0$. There exist $p'_n\in[0,p_n]$ such that
$$\sum p'_n=\nf\qquad\sum p'_na_n<\nf.$$
\end{lemma}
\par\noindent{\bf Proof } We may and shall assume that $\sum p_na_n=\nf$. There are only finitely many terms $a_n>1$. Replace these by $a'_n=1$. This has no effect on the convergence of $\sum p'_na_n$. Similarly we replace $a_n\in[1/2^k,2/2^k)$ for $k=1,2,\ldots$ by $a'_n=1/2^k$. Then $\sum p'_na_n$ converges if and only if $\sum p'_na'_n$ converges. Let $I_k$ be the set of indices $n$ for which $a'_n=1/2^k$, and $P_k$ the sum of $p_n$ over $n\in I_k$. Then $\sum P_k/2^k=\sum p_na'_n=\nf$. Hence there are infinitely many terms $P_k>1$. Set $p'_n=p_n/P_k$ for $n\in I_k$ if $P_k>1$. Then $P'_k=\min(1,P_k)$ and $\sum P'_k=\nf$. Also $\sum p'_na'_n=\sum P'_k/2^k\le\sum1/2^k=2$. \qed

\section{An application}

A classic result of Gnedenko states that the partial maxima of an iid sequence of positive random variables may be scaled to converge to 1 in probability, $M_n/a(n)\top1$, if the tail of the df varies rapidly, see~\cite{G43}. Resnick and Tomkins in~\cite{RT73} show that for any $c>1$ there exist dfs such that $\limsup M_n/a(n)=c$ almost surely. In~\cite{dHH72} it is shown that $\limsup M_n/a(n)=\nf$ a.s.\ is also possible.

\begin{theorem}\label{ta}
Let $G$ be a continuous strictly increasing df on $(0,\nf)$. One may choose $G$ such that the partial maxima $M_n$ satisfy $M_n/a(n)\top1$ with $1-G(a(n))=1/n$, and $\liminf M_n/a(n)=0$ almost surely.
\end{theorem}
\par\noindent{\bf Proof } Let $m_n=e^{s_n}$ be indices which increase so fast that $\s_n=s_n-s_{n-1}\to\nf$. Let $M_n$ denote the partial maxima from the sequence of independent observations $U_n$ from $G$. Set $m'_n=[\sqrt{m_nm_{n-1}}]$ and $x_n=a(m'_n)$. Define the events $B_n, E_n, A_n$ by: $B_n$ occurs if no observation $U_k$, $k\le m_n$, exceeds $x_n$, $E_n$ if $U_k>x_n$ holds for some $k\le m_{n-1}$ and $A_n$ if $U_k\le x_n$ for $m_{n-1}<k\le m_n$. Then $B_n=A_n\sm E_n$ and if one can show that $\sum\pbb B_n=\nf$, $\pbb E_n\to0$ and $x_n/a(m_n)\to0$ then almost surely $M_{m_n}/a(m_n)<x_n/a(m_n)$ infinitely often and hence $\liminf M_n/a(n)=0$ almost surely. 

First observe $\pbb E_n\le m_{n-1}\pbb\{U>x_n\}=m_{n-1}/m'_n\sim e^{-\s_n/2}\to0$. Now observe
$\pbb B_n=G(x_n)^{m_n}=e^{-\p_n}$ where $\p_n\sim m_n(1-G(x_n))$ since $m_n(1-G(x_n))^2\to0$. If $\p_n\sim\log\sqrt n$ then $\sum\pbb B_n=\nf$. Hence we choose $s_n=2n\log\log\sqrt n$. Then $\s_n=2\log\log\sqrt n+o(1)$ and $\p_n\sim m_n/m'_n=e^{\s_n/2}\sim\log\sqrt n$. We still have to choose $G$ such that $x_n/a(m_n)\to0$. Let $\q\in(0,1)$ and for $s>10$ write
\bql{qa}
t=T_0(s)=s/(\log\log s)^\q\qquad 1-G(e^t)=e^{-s};\qquad t_n=T_0(s_n)\qquad x_n=e^{t'_n}.
\eql
Then $t'_n=T_0(s_n-\s_n/2+o(1))$ and $t_n-t'_n\to\nf$. Hence $x_n/a(m_n)=e^{t'_n-t_n}\to0$. \qed

\end{document}